\documentclass[12pt, makeidx]{amsart}
\usepackage{amsthm}
\usepackage{amsfonts}
\usepackage{amssymb}
\newtheorem{thm}{Theorem}[]

\newtheorem{cor}[thm]{Corollary}
\newtheorem{lem}[thm]{Lemma}

\newtheorem{problem}[thm]{Problem}

\theoremstyle{definition}

\newtheorem{algorithm}[thm]{Algorithm}

\newtheorem{notation}[thm]{Notation}

\theoremstyle{remark}

\newtheorem{eg}[thm]{Example}

\numberwithin{figure}{section}
\newcommand{\ab}{{\rm ab}}
\newcommand{\Aut}{{\rm Aut(F)}}
\newcommand{\B}{\mathcal{B}}

\newcommand{\E}{{\mathcal E}}
\newcommand{\F}{{\rm F}}
\newcommand{\M}{{\rm M}}

\newcommand{\Integer}{\mathbb {Z}}

\newcommand{\Nat}{\mathbb {N}}

\newcommand{\abs}[1]{\left\vert#1\right\vert}

\newcommand{\normalclosure}[1]{\ll \!\! #1 \!\!\gg^{\rm F}}

\begin{document}

\title[Algorithmic constructions and primitive elements]{Algorithmic constructions and primitive elements in the free group of rank 2}


\author{Adam Piggott}

\address{School of Mathematics and Applied Statistics, University of Wollongong, NSW 2522, Australia}

\begin{abstract}
The centrepiece of this paper is a normal form for primitive
elements which facilitates the use of induction arguments to prove
properties of primitive elements. The normal form arises from an
elementary algorithm for constructing a primitive element $p$ in
$\F(x, y)$ with a given exponent sum pair $(X, Y)$, if such an
element $p$ exists. Several results concerning the primitive
elements of $\F(x, y)$ are recast as applications of the
algorithm and the normal form.
\end{abstract}

\maketitle

\section{Introduction}

Let $\F = \F(x, y)$ denote the free group on two generators $x$ and
$y$ and let $\F_{\ab} = \F_{\ab}(x, y)$ denote the free abelian
group on two generators $x$ and $y$. For an element $w \in \F$, the
\emph{exponent sum pair} is the ordered pair of integers $(X, Y)$
such that the exponent sum of $x$ in $w$ is $X$ and the exponent sum
of $y$ in $w$ is $Y$. Clearly, conjugate elements of $\F$ have the
same exponent sum pair. In the present paper, functions are written
to act on the right. Denote some specific automorphisms of $\F$ as
follows, where $v$ is
an element of $\F$:\\
\begin{center}
\begin{tabular}{cclccclccclcccl}
   & $\alpha_x$ &  &  & & $\alpha_y$ &  &  &  & $\beta$ &  &  &  & $\iota_v$ &   \\
  $x$ & $\mapsto$ & $x^{-1}$  &  & $x$ & $\mapsto$ & $x$      &  & $x$ & $\mapsto$ & $y$ &  & $x$ & $\mapsto$ & $v^{-1} x v$ \\
  $y$ & $\mapsto$ & $y$       &  & $y$ & $\mapsto$ & $y^{-1}$ &  & $y$ & $\mapsto$ & $x$ &  & $y$ & $\mapsto$ & $v^{-1} y v$ \\
\end{tabular}
\end{center}
An automorphism $\phi \in {\Aut}$ is said to be \emph{basic} if
there exists $n \in \Nat$ such that either $\phi$ is defined by $x
\mapsto xy^{n}$ and $y \mapsto xy^{n + 1}$ or $\phi$ is defined by
$x \mapsto xy^{n+1}$ and $y \mapsto xy^{n}$.  Let $\B$ denote the
set of basic automorphisms. Let $\Psi:\F \to \F$ be the map such
that $w\Psi = w^{-1} \alpha_x \alpha_y$ for each element $w \in \F$.

An element $w \in \F$ is said to be \emph{primitive} if it is the
image of $x$ under some automorphism $\theta_w \in \Aut$. Much is
known about the structure of primitive elements of $\F$. For
example, it was shown by Cohen, Metzler and
Zimmermann\cite{WhatDoesABasisLookLike} that, other than the
conjugacy class containing $x$ and the conjugacy class containing
$y$, each conjugacy class of primitive elements in $\F$ contains an
element of the form
$$x y^{m_1} x y^{m_s} \dots x y^{m_s} \eqno (\ast)$$ where $s \geq 0$ and
$m_i \in \{n, n+1\}$ for some $n \in \Nat$, or contains an element
obtained from an element of the form $(\ast)$ by application of some
combination of $\alpha_x$, $\alpha_y$ and $\beta$; we shall refer to
this fact as the \emph{(first) normal form property (for primitive
elements)}. An element $p \in \F$ is said to be a \emph{palindrome}
if $p \; \Psi = p$ (that is, ``$p$ reads the same forwards and
backwards''). It has recently been shown that each conjugacy class
of primitive elements in $\F$ contains an element $a$ such that
either $xay^{-1}$ is a palindrome or $y ax^{-1}$ is a
palindrome\cite[Theorem on p.613]{Helling}, and further that each
primitive element in $\F$ is the product of at most two
palindromes\cite[Lemma 1.6]{PrimitiveWidth}.

A theme of the present paper is the analysis of exponent sum pairs
to inform about primitive elements.  Such methods have been applied
since the seminal work of Nielsen in the early 20th century.

It is observed in \cite{WhatDoesABasisLookLike} that an elementary
algorithm for determining whether or not a particular element $w \in
\F$ is primitive follows from the normal form property. The
algorithm is modelled on the second of two proofs of the normal form
theorem, and provides evidence of the fundamental role that basic
automorphisms play in understanding primitive elements in $\F$.
Taking inspiration from Cohen, Metzler and Zimmermann's insight,
this paper records an algorithm which was developed from the
algorithm in \cite{WhatDoesABasisLookLike} and which solves the
following problem:

\begin{problem}\label{TheQuestion}
For relatively prime integers $X, Y \in \Integer$, write down a
primitive element $p \in \F$ with exponent sum pair $(X, Y)$.
\end{problem}

The utility of the above result is framed by the following two
well-known results.

\begin{lem}\label{DivisorsAreBad}
If $(X, Y)$ is the exponent sum pair of a primitive element $w$ in
$\F$, then $X$ and $Y$ are relatively prime.
\end{lem}

\begin{proof}
The element $w$ projects to $w_{\ab} = x^X y^Y$ in $\F_{\ab}$. Since
$w$ is primitive in $\F$, $w_{\ab}$ is primitive in $\F_{\ab}$ and
the result follows from the well-known analogous result in
$\F_{\ab}$.
\end{proof}

\begin{lem}[Nielsen, see {\cite[pp. 166-169]{MKS}})]\label{NielsenResult}
Each conjugacy class of primitive elements is determined uniquely by
the corresponding exponent sum pair.
\end{lem}

A proof of Lemma \ref{NielsenResult} is provided in
$\S$\ref{NormalClosureSection}.

Combining the solution to Problem \ref{TheQuestion} and Lemma
\ref{DivisorsAreBad} immediately yields the following.
\begin{thm}\label{ExponentSumPairsRelativelyPrime} There exists a primitive element $p$ in
$\F$
with exponent sum pair $(X, Y)$ if and only if $X$ and $Y$ are relatively prime integers.
\end{thm}

Combined with a simple observation and Lemma \ref{NielsenResult},
the solution to Problem \ref{TheQuestion} suggests another type of
normal form for primitive elements --- one which, for each primitive
$p \in \F$, describes an automorphism with the property that $y
\mapsto p$.


\begin{thm} [Second normal form for primitive elements]\label{InductiveStructure} For each \\primitive element $p
\in \F$ with exponent sum pair $(X, Y)$, there exist unique $\epsilon, \gamma, \delta \in \{0, 1\}$,
a unique minimal length $v \in \F$ such that the following conditions hold:
\begin{enumerate}
\item if $\abs{X} + \abs{Y} = 1$, then $p = y \; {\alpha_y}^\delta \beta^\epsilon \iota_v$;
\item if $\abs{X} + \abs{Y} = 2$, then $p = y \; \phi_\ast {\alpha_x}^\gamma {\alpha_y}^\delta \beta^\epsilon \iota_v$
        where $\phi_\ast \in \B$ is the basic automorphism such that $x \mapsto xy^2$ and $y \mapsto xy$;
\item if $\abs{X} + \abs{Y} > 2$, then there exist exactly two sequences of basic automorphisms $\phi_0, \phi_1, \dots,
\phi_s \in \B$ such that
$$p = y \; \phi_s \phi_{s-1} \dots \phi_0 {\alpha_x}^\gamma {\alpha_y}^\delta \beta^\epsilon \iota_v.$$
\end{enumerate}
Further, the values $\epsilon, \gamma, \delta, s$, the element $v \in \F$ and the basic
automorphisms $\phi_0, \phi_1, \dots, \phi_s \in \B$ may be found in
time proportional to $\log_2 \abs{p}$, where $\abs{p}$ denotes the
word-length of $p$.
\end{thm}

The second normal form confirms the importance of basic
automorphisms and offers a useful new perspective on the primitive
elements in $\F$. In particular, the second normal form facilitates
the use of inductive arguments (inducting on $s$) when proving
properties of primitive elements. Although such arguments are rarely
elegant, they are simple to implement. For example,
inductive arguments may be used to reprove the results from \cite{Helling} and
\cite{PrimitiveWidth} mentioned above.  The use of an inductive argument and the second
normal form provides common ground between these results, which at first sight appear
to be unrelated.

Let $r \in {\rm F}$. An algorithm for finding cyclically reduced
primitive $p$ such that $r$ is contained in the normal closure of
$p$ follows immediately from Algorithm \ref{TheAlgorithm} and the
following result.

\begin{thm}\label{AlgorithmToFindPforR}
Let $r \in \F$ and let $(A, B)$ be the exponent sum pair of $r$.  If
$(A, B) = (0, 0)$, then $r \in \normalclosure{p}$ for each primitive
element $p$ in $\F$.  If $(A, B) \neq (0, 0)$, then $r \in
\normalclosure{p}$ for primitive $p$ in $\F$ if and only if the
exponent sum pair of $p$ is
$\pm\frac{1}{d}(A, B)$ for $d$ the greatest common divisor of $A$ and $B$.
\end{thm}

The structure of this paper is as follows: in
$\S$\ref{AlgorithmSection} a solution to Problem \ref{TheQuestion}
is described and the second normal form theorem proved; in
$\S$\ref{ApplicationsOfSecondNormalForm} some applications of the
second normal form are detailed, including new proofs of those
results in \cite{PrimitiveWidth} \cite{Helling} described above; in
$\S$\ref{NormalClosureSection} Theorem \ref{AlgorithmToFindPforR} is
proved.

\section{A Solution to Problem \ref{TheQuestion}}\label{AlgorithmSection}

Let $\E$ denote the map from $\F$ to the set of ordered pairs of
integers, which maps $w \in \F$ to the exponent sum pair of $w$. Let
$\M$ denote the map ${\rm \Aut} \mapsto {\rm GL}(2, \Integer)$ such
that, for each automorphism $\theta \in {\rm Aut(F)}$,
$$\theta M = \left( \begin{array}{cc} X_1 & X_2 \\   Y_1 & Y_2 \\
\end{array}\right),$$ where $(X_1, Y_1)$ is the exponent sum pair of $x \theta$ and
$(X_2, Y_2)$ is the exponent sum pair of $y \theta$.  Let ${\rm
GL}(2, \Integer)$ act on the set of ordered pairs of integers by
matrix post-multiplication (where for this purpose an ordered pair
of integers is regarded as a $1 \times 2$ matrix of integers). It is
easily verified that $w \theta \E = (w \E) . (\theta \M)$ for each
automorphism $\theta \in \Aut$ and each $w \in \F$.






\begin{notation}
For integers $X, Y$ with $X \neq 0$, write $Y \!\!\! \mod X$ for the
unique integer $r$ such that $0 \leq r < \abs{X}$ and there exists
$q \in \Integer$ such that $Y = qX + r$.
\end{notation}

\begin{algorithm}\label{TheAlgorithm}
Let $(X, Y)$ be an ordered pair of relatively prime integers such
that $1 \leq X < Y$. Define $(X_0, Y_0) := (X, Y)$. Inductively, for
$i = 0, 1, 2, \dots$, proceed as follows:
\begin{itemize}
\item if $X_i = 1$ then terminate the inductive process;
\item if $X_i \geq 2$ and $Y_i \!\!\! \mod X_i \leq X_i - Y_i \!\!\! \mod X_i$ then define
    \begin{eqnarray}
    X_{i+1} & := & Y_i \!\!\!\!\! \mod{X_i} \nonumber \\
    Y_{i+1} & := & X_i - Y_i \!\!\!\!\! \mod{X_i} \nonumber \\
    n_i     & := & \max\{n \in \Nat \; | \; n X_i < Y_i \} \nonumber \\
    \phi_i & \in & \B \hbox{ such that } x \mapsto x y^{n_i+1} \hbox{ and } y  \mapsto x y^{n_i}; \nonumber
    \end{eqnarray}
\item if $X_i \geq 2$ and $Y_i \!\!\! \mod X_i > X_i - Y_i \!\!\! \mod X_i$ then define
    \begin{eqnarray}
    X_{i+1} & := & X_i - Y_i \!\!\!\!\! \mod{X_i} \nonumber \\
    Y_{i+1} & := & Y_i \!\!\!\!\! \mod{X_i} \nonumber \\
    n_i     & := & \max\{n \in \Nat \; | \; n X_i < Y_i \} \nonumber \\
    \phi_i & \in & \B \hbox{ such that } x \mapsto x y^{n_i} \hbox{ and } y  \mapsto x y^{n_i+1}. \nonumber
    \end{eqnarray}
\end{itemize}
\noindent It is clear that this inductive process terminates after
at most $\log_2 X$ iterations.  Let $s$ be the final value of $i$
considered.  The element $p := (xy^{Y_s})\phi_{s-1} \phi_{s-2} \dots
\phi_0$ is a primitive element in $\F$ with exponent sum pair $(X,
Y)$.
\end{algorithm}

\begin{proof}
It suffices to prove the following two claims, the first of which
confirms that the algorithm is well-defined and the second that it
achieves it goal.
\begin{enumerate}
\item [(A)] for each integer $i = 0, \dots, s-2$, if $(X_i, Y_i)$
is an ordered pair of relatively prime integers then $(X_{i+1}, Y_{i+1})$
is also an ordered pair of relatively prime integers;
\item [(B)] for each integer $i = s-1, \dots, 1$, if $p_{i+1}$ is a primitive element with
exponent sum pair $(X_{i+1}, Y_{i+1})$, then $(p_{i+1})\phi_i$ is a
primitive element with exponent sum pair $(X_i, Y_i)$.
\end{enumerate}    Both claims are proved by
inductive arguments.  The inductive steps are shown below.

Let $i$ be an integer such that $0 \leq i < s-1$ and assume that
$(X_i, Y_i)$ is an ordered pair of relatively prime integers.
Suppose that $d$ is a positive integer such that $d$ divides both
$X_{i+1}$ and $Y_{i+1}$. By definition, $d$ divides both $Y_i \!\!\!
\mod X_i$ and $(X_i - Y_i \!\!\! \mod X_i)$.  Since $X_i = (X_i -
Y_i \!\!\! \mod X_i) + Y_i \!\!\! \mod X_i$, it follows that $d$
divides $X_i$. Since $Y_i = n_i X_i + Y_i \!\!\! \mod X_i$, it
follows that $d$ divides $Y_i$. Now, $d$ divides both $X_i$ and
$Y_i$ and $X_i, Y_i$ relatively prime implies that $d = 1$, hence
$(X_{i+1}, Y_{i+1})$ is also an ordered pair of relatively prime
integers.

The inductive step in the proof of Claim (B) is easily verified by
calculation as follows. In the case that $Y_i \!\!\! \mod X_i \leq
X_i - Y_i \!\!\! \mod X_i$, then
\begin{eqnarray*}
(p_{i+1})\phi_i \E & = & (p_{i+1} \E) . (\phi_i \M) \\
  & = & (X_{i+1}, Y_{i+1}) \left(%
\begin{array}{cc}
  1 & n_i \\
  1 & n_i+1 \\
\end{array}%
\right) \\
 & =&  (X_{i+1} + Y_{i+1}, n_i(X_{i+1} + Y_{i+1}) + Y_{i+1}) \\
 & =& (X_i, Y_i). \\
\end{eqnarray*}  In the case that $Y_i \!\!\! \mod X_i > X_i - Y_i \!\!\! \mod
X_i$, then
\begin{eqnarray*}
(p_{i+1})\phi_i \E & = & (p_{i+1} \E) . (\phi_i \M) \\
 & = & (X_{i+1}, Y_{i+1}) \left(%
\begin{array}{cc}
  1 & n_i+1 \\
  1 & n_i \\
\end{array}%
\right) \\
 & =&  (X_{i+1} + Y_{i+1}, n_i(X_{i+1} + Y_{i+1}) + X_{i+1}) \\
 & =& (X_i, Y_i). \\
\end{eqnarray*}
\end{proof}


Algorithm \ref{TheAlgorithm} is easily extended to all relatively
prime ordered pairs of integers $(X, Y)$, and hence a solution to
Problem \ref{TheQuestion}, by the following observations:
\begin{enumerate}
\item it follows from Lemma \ref{DivisorsAreBad} that there is no primitive element in $F$ with exponent sum pair $(0, 0)$;
\item \label{SignNoProblem} it follows from the properties of the automorphisms $\alpha_x, \alpha_y, \beta \in {\rm Aut(F)}$,
that there exists a primitive element in $F$ with exponent sum pair
$(X, Y)$ if and only if there exists a primitive element in $F$ with
exponent sum pair $(\min\{\abs{X}, \abs{Y}\}, \max\{\abs{X},
\abs{Y}\})$.
\end{enumerate}

\begin{eg}
Find a primitive element with exponent sum pair $(34, -27)$.
\vskip10pt \noindent Define $X_0 := 27$, $Y_0 := 34$.

    Since $7 = Y_0 \!\!\! \mod X_0 < X_0 - Y_0 \!\!\! \mod X_0 = 20$, define $X_1 := 7$, $Y_1 := 20$, $n_0 := 1$ and
    $\phi_0 \in {\rm Aut(F)}$ such that $x \mapsto xy^2$ and $y \mapsto xy$.

    Since $6 = Y_1 \!\!\! \mod X_1 > X_1 - Y_1 \!\!\! \mod X_1 = 1$, define $X_2 := 1$, $Y_2 := 6$, $n_1 := 2$ and
    $\phi_1 \in {\rm Aut(F)}$ such that $x \mapsto xy^2$ and $y \mapsto xy^3$.

    Then $ (x y^6) \phi_1 \phi_0 = \bigl(xy^2 (xy^3)^6\bigr) \phi_0  = xy^2 (xy)^2 \bigl( xy^2 (xy)^3\bigr)^6$ is a primitive element with exponent sum pair
    $(27, 34)$, and $$(x y^6)\phi_1 \phi_0 \alpha_x \beta  = y^{-1}x^2 (y^{-1}x)^2 \bigl( y^{-1} x^2 (y^{-1} x)^3\bigr)^6$$ is a primitive element
    with exponent sum pair $(34, -27)$.
\end{eg}

To prove the second normal form theorem, it is convenient to use the
following lemma.

\begin{lem}\label{UniqueBasicSequence}
Let $(X, Y)$ be an ordered pair of relatively prime natural numbers such that $1 \leq X < Y$.
There exists a unique sequence of ordered pairs $(1, Y_s) = (X_s, Y_s), (X_{s-1}, Y_{s-1}), \dots (X_0, Y_0) = (X, Y)$
such that, for each $i = s-1, \dots, 1, 0$:
$$(X_{i+1}, Y_{i+1}). (\phi_i M) = (X_i, Y_i),$$ for some basic automorphism $\phi_i \in \B$.
\end{lem}

\begin{proof}
It is easily verified that for each basic automorphism $\phi \in \B$ and each ordered pair $(U, V)$,
either $(U, V). (\phi_i M) = (U+V, n(U+V)+U)$ or $(U, V). (\phi_i M) = (U+V, n(U+V)+V)$.  In either case,
if $(X_{i+1}, Y_{i+1}). (\phi_i M) = (X_i, Y_i)$, then $X_{i+1} = \min\{Y_i \!\!\! \mod X_i, X_i - Y_i \!\!\! \mod X_i\}$ and
$Y_{i+1} = \max\{Y_i \!\!\! \mod X_i, X_i - Y_i \!\!\! \mod X_i\}$.  Thus the sequence $(1, Y_s) = (X_s, Y_s),
(X_{s-1}, Y_{s-1}), \dots (X_0, Y_0) = (X, Y)$ determined in Algorithm $\ref{TheAlgorithm}$ is the unique sequence
with the desired properties.
\end{proof}

The second normal form is proved by collating some of the results
obtained above.

\begin{proof}[Proof of the second normal form theorem]
By the properties of $\alpha_x$, $\alpha_y$, $\beta$ and the set of
inner automorphisms, it suffices to consider cyclically reduced
primitive elements $p \in \F$ with exponent sum pairs $(X, Y)$ such
that $0 \leq X \leq Y$. In the case that $\abs{X} + \abs{Y} = 1$, Lemma
\ref{DivisorsAreBad} implies that $X = 0$, $Y = 1$ and $p$ is a cyclic
permutation of $(y) {\alpha_x}^0 {\alpha_y}^0 \beta^0 \iota_1$.  In
the case that $\abs{X} + \abs{Y} = 1$, then $X = 1$ and $Y = 1$ and the result is clear.
In the case that $X = 1$ but $\abs{X} + \abs{Y} > 2$, define $\phi_0 \in \B$ such that
$y \mapsto xy^Y$ and either $x \mapsto xy^{Y-1}$ or $x \mapsto xy^{Y+1}$, then $p$ is a cyclic permutation of
$(y) \phi_0 {\alpha_x}^0 {\alpha_y}^0 \beta^0 \iota_1$. In the case
that $X \geq 2$, Lemma \ref{DivisorsAreBad} implies that $X < Y$.  It follows from the algorithm
and Lemma \ref{UniqueBasicSequence} that there is a unique sequence of exponent sum pairs
$(1, Y_s) = (X_s, Y_s), (X_{s-1}, Y_{s-1}), \dots (X_0, Y_0) = (X, Y)$ such that, for each $i = s-1, \dots, 1, 0$,
$(X_{i+1}, Y_{i+1}). (\phi_i M) = (X_i, Y_i),$ for some basic automorphism $\phi_i \in \B$.  If $Y_s = 1$, then $\phi_{s-1}$
is such that $x \mapsto xy^2$ and $y \mapsto xy$, but $\phi_{s-2}$ may be defined such that $x \mapsto xy^{n_i}$ and
$y \mapsto xy^{n_i+1}$ or $x \mapsto xy^{n_i+1}$ and
$y \mapsto xy^{n_i}$; it is clear that the remaining basic automorphisms are uniquely determined by the sequence of
exponent sum pairs.  If $Y_s > 1$, then $\phi_{s-1}$
may be defined such that $y \mapsto xy^Y_s$ and
$x \mapsto xy^{Y_s-1}$ or $x \mapsto xy^{Y_s+1}$; it is clear that the remaining basic
automorphisms are uniquely determined by the sequence of exponent sum pairs.
\end{proof}

\section{Some Applications of the Second Normal Form}\label{ApplicationsOfSecondNormalForm}

In this section, some applications of the second normal form and are
described. It is convenient to first record the following lemma, the
proof of which is trivial.

\begin{lem}\label{PalindromeLemma}
Let $\phi$ be a basic automorphism.  If $w \in \F$ is a palindrome
in which only positive exponents appear, then $(w)\phi = x v$ for
some palindrome $v \in \F$ in which only positive exponents appear.
\end{lem}

It is now possible to reprove the result from \cite{PrimitiveWidth}
mentioned in the introduction, using an induction technique based on
the second normal form.

\begin{thm}[Shpilrain, Bardakov and Tolstykh \cite{PrimitiveWidth}]\label{SBTTheorem}
Each primitive element $p \in \F$ is either a palindrome, or is
the product of two palindromes.
\end{thm}

\begin{proof}
Let $\epsilon, \gamma, \delta, v, s$, and $\phi_i$ (for $i = 0, 1,
\dots s$) be as in the statement of Corollary
\ref{InductiveStructure}.  If $s = -1$, then $p \in \{x^{\pm 1},
y^{\pm 1}\}$ and $p$ is a palindrome.  If $s = 0$, then $p \in
\{x^{\pm 1} y^z \; | \; z \neq 0\} \cup \{x^z y^{\pm 1} \; | \; z
\neq 0\}$ and $p$ is the product of two palindromes.  Assume the
result holds for each primitive element where $s = k$, for some $k
\geq 0$. Consider the case that $s = k+1$.  By the inductive
hypothesis, $(y)\phi_s \dots \phi_1$ is a either a palindrome or a
product of two palindromes.  In the former case, Lemma
\ref{PalindromeLemma} informs that $(y) \phi_s \dots \phi_0 = x v$
for some palindrome $v \in \F$; in the latter case, say $(y)\phi_s
\dots \phi_1 = v_1 v_2$ for palindromes $v_1, v_2 \in \F$, Lemma
\ref{PalindromeLemma} informs that $(y) \phi_s \dots \phi_0 = x v_3
x v_4 = (x v_3 x) v_4$ for some palindromes $v_3, v_4 \in \F$. Hence
in either case, $(y)\phi_s \dots \phi_0$ is the product of two
palindromes (and possibly also a palindrome itself). It is clear
that application of $\alpha_x$, $\alpha_y$, $\beta$ and inner
automorphisms preserve the property of being a palindrome or being a
product of two palindromes, hence $p$ has the required property.
\end{proof}

The author is grateful to Peter Nickolas for pointing out the
following corollary to Theorem \ref{SBTTheorem}.

\begin{cor}
Let $p$ be a primitive element in $\F$.  One of the following
two statements holds:
\begin{enumerate}
\item \label{CASE1} $p = z^{-1} w z$ for some $z \in \F$ and some palindrome $w \in \F$;
\item \label{CASE2} $p = z^{-1} a w z$ for some $z \in \F$, some $a \in \{x, x^{-1}, y, y^{-1}\}$ and some palindrome $w \in \F$.
\end{enumerate}
\end{cor}

\begin{proof}
Let $p$ be a primitive element in $\F$.  If $p$ is a palindrome,
there is nothing to prove, so we may assume that $p$ is not a
palindrome.  By the Theorem, there exist palindromes $w_1, w_2 \in
\F$ such that $p = w_1 w_2$. Consider first the case that $w_1$ has
even length, say $w_1 = v (v \Psi)$ for some $v \in \F$.  Then $v
((v \Psi) w_2 v) v^{-1} = v (v \Psi) w_2 = p$; hence Case
(\ref{CASE1}) holds with $z = v^{-1}$. Next, consider the case that
$w_1$ has odd length, say $w_1 = v a (v \Psi)$ for some $v \in \F$
and some $a \in \{x, x^{-1}, y, y^{-1}\}$.  Then $v (a (v \Psi) w_2
v) v^{-1} = v a (v \Psi) w_2 = p$; hence Case (\ref{CASE2}) holds
with $z = v^{-1}$.
\end{proof}

We may use a similar strategy to reprove the result from
\cite{Helling} mentioned in the introduction.

\begin{lem}[Helling \cite{Helling}]
For each primitive element $p \in \F$, there exists a
palindrome $v \in \F$ and an element $z \in \F$ such that
either $p = z y^{-1} v x z^{-1}$ or $p = z x^{-1} v y z^{-1}$.
\end{lem}

\begin{proof}
Define a subset $\mathcal{Q}$ of $\F$ as follows:

$$\mathcal{Q} := \{w \in \F \; | \; \exists z \in \F \hbox{ and a palindrome } v \in \F \hbox{ such that }$$
$$ w = z y^{-1} v x z^{-1} \hbox{ or } w = z x^{-1} v y z^{-1}\}.$$
It is clear from the definition that $\mathcal{Q}$ is closed under
the action of inner automorphisms.

Let $p$ be such that $p = z y^{-1} v x z^{-1}$ for some $z \in \F$ and for some palindrome $v \in \F$.  Then
\begin{eqnarray}
(p)\alpha_x     & = &   (z \alpha_x) y^{-1} (v \alpha_x) x^{-1} {(z^{-1} \alpha_x)} \nonumber \\
                & = &   (z \alpha_x) (x x^{-1}) y^{-1} (v \alpha_x) (y^{-1} y) x^{-1} {(z^{-1} \alpha_x)} \nonumber \\
                & = &   ((z \alpha_x) x) x^{-1} (y^{-1} (v \alpha_x) y^{-1}) y (x^{-1} {(z^{-1} \alpha_x)}) \nonumber \\
                & = &   z' x^{-1} u y (z')^{-1}, \nonumber
\end{eqnarray} where $u = y^{-1} (v \alpha_x) y^{-1}$ is a palindrome and $z' = (z \alpha_x) x$;
\begin{eqnarray}
(p)\alpha_y     & = &   (z \alpha_y) y (v \alpha_y) x {(z^{-1} \alpha_y)} \nonumber \\
                & = &   (z \alpha_y) y (x^{-1} x) (v \alpha_y) x (y y^{-1}){(z^{-1} \alpha_y)} \nonumber \\
                & = &   ((z \alpha_y) y) x^{-1} (x (v \alpha_y) x) y (y^{-1}{(z^{-1} \alpha_y})) \nonumber \\
                & = &   z' x^{-1} u y (z')^{-1} \nonumber
\end{eqnarray} where $u = x (v \alpha_y) x$ is a palindrome and $z' = (z \alpha_y) y$; and
\begin{eqnarray}
(p)\beta        & = &   (z \beta) x^{-1} (v \beta) y {(z^{-1} \beta)} \nonumber \\
                & = &   z' x^{-1} u y (z')^{-1} \nonumber
\end{eqnarray} where $u = v \beta$ is a palindrome and $z' = z \beta$.  A similar treatment shows that if $p$
is such that $p = z x^{-1} v y z^{-1}$, then $p \alpha_x$, $p
\alpha_y$ and $p \beta$ are of the form $z' x^{-1} u y (z')^{-1}$ or
$z' y^{-1} u x (z')^{-1}$, for some palindrome $u \in \F$ and for
some $z' \in \F$.  Hence $\mathcal{Q}$ is also closed under the
action of $\alpha_x, \alpha_y$, and $\beta$.

To complete the proof it suffices to show that $(y)\phi_s \dots
\phi_0 \in \mathcal{Q}$ for each list of basic automorphisms
$\phi_0, \phi_1, \dots, \phi_s$.  Let $s \geq 0$ and let $\phi_0,
\phi_1, \dots, \phi_s$ be a list of basic automorphisms. It is clear
that $(y) \phi_s \in \mathcal{Q}$. Assume that, for some $i$ such
that $0 < i \leq s$, $(y)\phi_s \dots \phi_i \in \mathcal{Q}$.
Suppose first that
 $(y)\phi_s \dots \phi_i = z x^{-1} v y z^{-1}$ for some palindrome $v \in \F$ and for
some $z \in \F$ and $\phi_{i-1}$ is defined by $x \mapsto
xy^{n_i+1}$ and $y \mapsto xy^{n_i}$.  Then
\begin{eqnarray}
(y) \phi_s \dots \phi_i \phi_{i-1}(y) & = & (z x^{-1} v y z^{-1})\phi_{i-1} \nonumber \\
    & = & (z \phi_{i-1}) y^{-(n_i+1)}x^{-1} x u x y^{n_i} {(z^{-1} \phi_{i-1})}) \nonumber \\
    & &     \hbox{(for some palindrome u}\!\in \F \hbox{, by Lemma \ref{PalindromeLemma})}\nonumber \\
    & = & ((z \phi_{i-1}) y^{-n_i}) y^{-1} u x (y^{n_i} {(z^{-1} \phi_{i-1})}) \nonumber \\
    & = & z' y^{-1} u x (z')^{-1}, \nonumber
\end{eqnarray} where $z' = (z \phi_{i-1}) y^{-n_i}$.  Now suppose that
$\phi_{i-1}$ is defined by $x \mapsto xy^{n_i}$ and $y \mapsto
xy^{n_i+1}$.  Then
\begin{eqnarray}
(y) \phi_s \dots \phi_i \phi_{i-1} & = & (z x^{-1} v y z^{-1})\phi_{i-1} \nonumber \\
    & = & (z \phi_{i-1}) y^{-n_i}x^{-1} x u x y^{n_i+1} {(z^{-1} \phi_{i-1})}) \nonumber \\
    & &     \hbox{(for some palindrome u}\!\in \F \hbox{, by Lemma \ref{PalindromeLemma})}\nonumber \\
    & = & ((z \phi_{i-1}) y^{-n_i}) x^{-1} (x u x) y (y^{n_i} {(z^{-1} \phi_{i-1})}) \nonumber \\
    & = & z' x^{-1} (x u x) y (z')^{-1}, \nonumber
\end{eqnarray} where $x u x$ is a palindrome and $z' = (z \phi_{i-1}) y^{-n_i}$.  The case that
$(y)\phi_s \dots \phi_i = z y^{-1} v x z^{-1}$ for some palindrome
$v \in \F$ and for some $z \in \F$, is verified similarly.
\end{proof}

\section{The Normal Closure of a Primitive Element}\label{NormalClosureSection}


Let $\{p, q\}$ be a basis for $\F$.  For each element $r \in \F$,
let $r_{p, q}$ denote the unique reduced word in $\{p^{\pm 1},
q^{\pm 1}\}$ such that $r_{p, q}$ is equal to $r$ in $F$. The
\emph{normal closure} of $p$ in $\F$, denoted $\normalclosure{p}$,
is defined to be
$$\normalclosure{p} := \{ w \in \F \; | \; w = \overset{s}{\underset{i = 1}{\prod}}{f_i}^{-1} p^{\epsilon_i}f_i, \hbox{ for some } s \geq 0, \epsilon_j \in {\pm 1}, f_j \in \F\}.$$

\begin{lem}\label{NormalClosureAndExponentSum}
For each $r \in \F$, $r \in \normalclosure{p}$ if and only if the
exponent sum of $q$ in $r_{p, q}$ is zero.
\end{lem}

\begin{proof}
Suppose that the exponent sum of $q$ in $r_{p, q}$ is zero. Then
$$r_{p, q} = p^{\alpha_1} q^{\beta_1} \dots p^{\alpha_s} q^{\beta_s},$$
for some $s \geq 0$, $\alpha_j$ non-zero (except perhaps
$\alpha_1$), $\beta_j$ non-zero (except perhaps $\beta_s$) such that
$\underset{j=1}{\overset{s}{\sum}} \beta_j = 0$.  Insertion of
trivial words yields
$$ r_{p, q} = \bigl(q^{-(\beta_1 + \dots + \beta_s)} p^{\alpha_1} q^{\beta_1 + \dots \beta_s}\bigr) \dots \bigl(q^{-\beta_{s-1}-\beta_s} p^{s-1}
q^{\beta_{s-1} + \beta_s}\bigr) \bigl(q^{-\beta_s}p^{\alpha_s}
q^{\beta_s}\bigr),$$ and $r \in \normalclosure{p}$. The opposite
direction of implication follows easily from the definition of
$\normalclosure{p}$.
\end{proof}

As an aside to ensure that the present paper is self-contained,
Lemma \ref{NormalClosureAndExponentSum} may be used to prove
Nielsen's result, Lemma \ref{NielsenResult}.
\begin{proof}[Proof of Lemma \ref{NielsenResult}]
Let $z \in \F$ be a primitive element with exponent sum pair
$(0, 1)$.  By Lemma \ref{NormalClosureAndExponentSum}, $z \in
\normalclosure{y}$ and $y \in \normalclosure{z}$.  That is,
$$z = \overset{s}{\underset{i = 1}{\prod}} {w_i}^{-1} y^{\epsilon_i} w_i, \hbox{ and }  y = \overset{t}{\underset{j = 1}{\prod}} {v_j}^{-1} z^{\delta_j} v_j,$$
for some $s, t \in \Nat$, some $w_i, v_j \in \F$ and some
$\epsilon_i, \delta_j \in \{\pm 1\}$ such that
$\overset{s}{\underset{i = 1}{\sum}}\epsilon_i = 1$ and
$\overset{t}{\underset{j = 1}{\sum}}\delta_j = 1$. Substitution
yields
\begin{eqnarray}
z & = & \overset{s}{\underset{i = 1}{\prod}} {w_i}^{-1} \bigl( \overset{t}{\underset{j = 1}{\prod}} {v_j}^{-1} z^{\delta_j} v_j \bigr) w_i \nonumber \\
  & = & {w_1}^{-1} {v_1}^{-1} z^{\delta_1} v_1  \dots {v_t}^{-1} z^{\delta_t} v_t w_1 \dots {w_s}^{-1} {v_1}^{-1} z^{\delta_1} v_1  \dots {v_t}^{-1} z^{\delta_t} v_t w_s.\nonumber
\end{eqnarray}
It follows that $v_j {v_{j+1}}^{-1} = 1$ for each $j = 1, \dots,
t-1$ and $w_i {w_{i+1}}^{-1} = 1$ for each
$i = 1, \dots, s-1$, hence we may assume that $s = 1$ (and $t = 1$) and $z$ is conjugate to $y$.  

More generally, let $p$ and $q$ be primitive elements with the same
exponent sum pair.  Since $p$ is primitive, there exists an
automorphism $\theta \in {\rm Aut(F)}$ such that $p \mapsto y$.  It
follows that the exponent sum pair of $q\theta$ is $(0, 1)$, hence
$q \theta$ is conjugate to $p \theta$ and $q$ is conjugate to $p$.
\end{proof}

Let $(X, Y)$ denote the exponent sum pair of $p$ and let $(U, V)$
denote the exponent sum pair of $q$.

\begin{lem}\label{DeterminantCondition}
It holds that $XV - YU = 1$.
\end{lem}

\begin{proof}
Since $\{p, q\}$ is a basis for $\F$, $\{x^Xy^Y, x^Uy^V\}$ is a
basis for $F_{\ab}$.  The result then follows from the well-known
analogous result for $F_{\ab}$.
\end{proof}

\begin{cor}\label{TheCor}
For each $r \in \F$, $r \in \normalclosure{p}$ if and only if
the exponent sum pair of $r$ is $(kX, kY)$ for some integer $k$.
\end{cor}

\begin{proof}
Let $r \in \F$, let $P$ denote the exponent sum of $p$ in $r_{p, q}$
and let $Q$ denote the exponent sum of $q$ in $r_{p, q}$.  Note that
the exponent sum pair of $r$ is given by $(PX + QU, PY+QV)$. It
follows easily from the definition of $\normalclosure{p}$, that $r
\in \normalclosure{p}$ implies the exponent sum pair of $r$ is $(PX,
PY)$.  Suppose that the exponent sum pair of $r$ is $(kX, kY)$ for
some integer $k$.  If $Q$ is non-zero, then
 $(PX + QU, PY+QV) = (kX, kY)$ implies that $U = \frac{(k-P)X}{Q}$ and $V = \frac{(k-P)Y}{Q}$ and $XV - YU = 0$,
contradicting Lemma \ref{DeterminantCondition}; hence $Q = 0$ and
Lemma \ref{NormalClosureAndExponentSum} implies that $r \in
\normalclosure{p}$.
\end{proof}

Theorem \ref{AlgorithmToFindPforR} follows immediately from
Corollary \ref{TheCor}.
\bibliographystyle{amsplain}
\bibliography{PrimitiveElementsBib}
\end{document}